\documentclass[12pt]{amsart}

\usepackage[T1]{fontenc}

\usepackage{verbatim}

\usepackage{tikz-cd}

\usepackage[backref,bookmarks=false]{hyperref}

\usepackage{color}

\numberwithin{equation}{section}

\newcommand{\ben}{\begin{enumerate}}
\newcommand{\een}{\end{enumerate}}

\newcommand{\bea}{\begin{eqnarray}}
\newcommand{\ba}{\begin{array}}
\newcommand{\bean}{\begin{eqnarray*}}
\newcommand{\ea}{\end{array}}
\newcommand{\eea}{\end{eqnarray}}
\newcommand{\eean}{\end{eqnarray*}}
\newcommand{\beq}{\begin{equation}}
\newcommand{\eeq}{\end{equation}}
\newcommand{\bthm}{\begin{thm}}
\newcommand{\ethm}{\end{thm}}
\newcommand{\blem}{\begin{lem}}
\newcommand{\elem}{\end{lem}}
\newcommand{\bprop}{\begin{prop}}
\newcommand{\eprop}{\end{prop}}
\newcommand{\bcor}{\begin{cor}}
\newcommand{\ecor}{\end{cor}}
\newcommand{\bdfn}{\begin{dfn}}
\newcommand{\edfn}{\end{dfn}}
\newcommand{\brem}{\begin{rem}}
\newcommand{\erem}{\end{rem}}
\newcommand{\bpf}{\begin{proof}}
\newcommand{\epf}{\end{proof}}
\newcommand{\bfact}{\begin{fact}}
\newcommand{\efact}{\end{fact}}
\newcommand{\bobs}{\begin{obs}}
\newcommand{\eobs}{\end{obs}}
\newcommand{\bexam}{\begin{exam}}
\newcommand{\eexam}{\end{exam}}
\newcommand{\bclaim}{\begin{claim}}
\newcommand{\eclaim}{\end{claim}}

\newtheorem{thm}{Theorem}[section]
\newtheorem{prop}[thm]{Proposition}
\newtheorem{lem}[thm]{Lemma}

\newtheorem{cor}[thm]{Corollary}
\newtheorem{dfn}[thm]{Definition}
\newtheorem{rem}[thm]{Remark}
\newtheorem{fact}[thm]{Fact}
\newtheorem{claim}[thm]{Claim}
\newtheorem{obs}[thm]{Observation}
\newtheorem{exam}[thm]{Example}

\newtheorem*{condition'}{Condition 2'}

 \newtheoremstyle{claimstyle}%
   {}
   {}
   {\normalfont}
   {}
   {\itshape}
   {.}
   { }
   {\thmnote{#3}}

\theoremstyle{claimstyle}

\renewcommand{\mod}{\operatorname{mod}}

\alph{enumii} \roman{enumiii}

\def\cA{\mathcal A}             \def\cB{\mathcal B}       \def\cC{\mathcal C}

\def\cU{\mathcal U}                       \newcommand{\J}{\mathcal{J}}
              \def\cS{\mathcal S}             
                    
 \def\cI{\mathcal I}

\def\N{{\mathbb N}}                \def\Z{{\mathbb Z}}      \def\R{{\mathbb R}}
\def\C{{\mathbb C}}                      \def\oc{{\widehat \C}}

\newcommand{\cbar}{\widehat{{\mathbb C}} }

                \def\b{\beta}             \def\d{\delta}
                          
                           \def\l{\lambda} 
              \def\om{\omega}           \def\Om{\Omega}

\def\ka{\kappa}

\newcommand{\ph}{\varphi}
\newcommand{\al}{\alpha}
\newcommand{\ga}{\gamma}

\def\1{1\!\!1}

\def\and{\text{ and }}

        \def\diam{\text{\rm {diam}}}

\def\Sing{S}

     \def\HD{\text{{\rm HD}}}

         \def\P{\text{{\rm P}}}

              \def\bu{\bigcup}
\def\({\bigl(}                \def\){\bigr)}
\def\lt{\left}                \def\rt{\right}

                        \def\^{\tilde}

\def\es{\emptyset}            \def\sms{\setminus}
\def\sbt{\subset}             \def\spt{\supset}

       \def\lra{\longrightarrow}

\def\sp{\medskip}                     
\def\ov{\overline}

\def\om{\omega}

\def\D{{\mathbb D}}

\def\${$ \displaystyle }

\newcommand{\cri}{{\mathcal C _f}}

\newcommand{\jul}{\mathcal J}
\newcommand{\fat}{\mathcal F}

\def\escf{\cI (f)}
\def\esc{\cI}

\def\hd{\text{{\rm HD}}}
\def\crex{\d}

\newcommand{\singf}{S(f)}
\newcommand{\sing}{S}

\def\alm{{\bf a}}
\def\blm{{\bf b}}
\def\phlm{{\bf  \ph}}
\def\Vlm{{\bf V}}
\def\Ulm{{\bf U}}

\def\II{\mathcal I}

\begin{document}

\title[The Hausdorff dimension of escaping sets]{The exact value of Hausdorff dimension \\ of \\ escaping sets \\ of \\ class $\cB$ meromorphic functions }


\author[\sc Volker MAYER]{\sc Volker Mayer}
\address{Volker Mayer, Universit\'e de Lille, D\'epartement de
  Math\'ematiques, UMR 8524 du CNRS, 59655 Villeneuve d'Ascq Cedex,
  France} \email{volker.mayer@univ-lille.fr \newline
  \hspace*{0.42cm} \it Web: \rm math.univ-lille1.fr/$\sim$mayer}


\author[\sc Mariusz Urba\'nski]{\sc Mariusz Urba\'nski}
\address{Mariusz Urba\'nski, 
Department of
  Mathematics, University of North Texas, Denton, TX 76203-1430, USA}
\email{urbanski@unt.edu \newline \hspace*{0.42cm} \it Web: \rm
  www.math.unt.edu/$\sim$urbanski}

\date{\today} 


\begin{abstract}  
We consider the subclass of class $\cB$ consisting of meromorphic functions $f:\C\to\cbar$ for which infinity is not an asymptotic value and whose all poles have orders uniformly bounded from above. This class  was introduced in \cite{BwKo2012} and the Hausdorff dimension $\HD(\escf)$ of the set $\escf$ of all points escaping to infinity under forward iteration of $f$ was estimated therein. In this paper we provide a closed formula for the exact value of $\HD(\escf)$ identifying it with the critical exponent of the natural series introduced in \cite{BwKo2012}. This exponent is very easy to calculate for many concrete functions. In particular, we construct a function from this class which is of infinite order and for which $\HD(\escf)=0$. 
 \end{abstract}

\maketitle



\section{Introduction} \label{s1}

We deal with the dynamics of meromorphic functions from class $\cB$. This class consists in all meromorphic functions $f:\C\to\cbar$ for which the 
set $\Sing (f)$ of finite asymptotical and critical values is bounded. For entire functions, class $\cB$ was introduced by Eremenko and Lyubich in \cite{EL92}.

Recalling standard definitions, the \emph{Fatou} set $\fat (f)$ of a function $f\in \cB$ consists in all points in $\C$ that admit a neighborhood on which all iterates of $f$ are defined and normal in the sense of Montel. Then, the \emph{Julia} set $\jul (f)$ of $f$ is the complement of the \emph{Fatou} set $\fat (f)$.

In this paper we focus on another dynamically significant set, namely the \emph{escaping} set
of $f$
\beq\label{5}
\escf =\{z\in \C \; : \; f^n(z)\to\infty \;\; \text{as}\;\; n\to\infty\}
\eeq
which is also  the intersection over all $R>0$ of the sets
\beq\label{6}
\esc_R(f) =\{z\in \C \; : \; \liminf_{n\to\infty} |f^n(z)|\geq R\}.
\eeq
If $f$ is a polynomial then $\escf$ is simply the basin of (immediate) attraction to infinity which is a central object in the study of the dynamics of polynomials.
If $f\in \cB$ is transcendental then $\escf\neq \emptyset$ and $\escf\subset \jul (f)$ (see \cite{EL92} and \cite{RS99a}) and we are interested in 
$$
\HD(\escf),
$$
where $\hd(E)$ denotes the \emph{Hausdorff dimension} of a set $E$ contained in a metric space. 


McMullen showed in \cite{McM87} that the Julia set of any exponential function $z\mapsto\l e^z$, $\l\neq 0$, has Hausdorff dimension equal to $2$
and that the Julia set of any sine function $z\mapsto\sin (\al z +\b )$, $\al\neq 0$, has even positive Lebesgue measure. His result has been extended by Bara\'nski 
\cite{B08} and Schubert \cite{Schubert} to all entire functions in
class $\cB$ that have finite order. 
As it is mentioned in \cite{BwKo2012}, their proof actually shows that 
$$\HD(\escf) =2$$
for all such functions and in fact for all meromorphic functions in the Eremenko-Lyubich class $\cB$ that have finite order and 
for which $\infty$ is an asymptotic value. Bergweiler, Karpi\'nska, and Stallard  showed in \cite{BwKaSta2009} that this result
even also holds for some functions $f\in \cB$ having infinite order. 
The situation is totally different when $\infty$ is not an asymptotical value as was shown by Bergweiler and Kotus in \cite{BwKo2012} who, improving former results for elliptic functions by Kotus and Urba\'nski (see \cite{KU03}, comp. \cite{KUbook}), obtained the following upper estimate.

\bthm[Bergweiler and Kotus]\label{BK theo}
Let $f\in \cB$ be a transcendental meromorphic function having finite order $\rho <\infty$, for which $\infty$ is not an asymptotic value and 
such that there exists $M \in \N$ such that the multiplicity of all poles, except possibly finitely many, is at most M. Then
\beq\label{7}
\HD(\escf) \leq \frac{2M\rho}{2+M\rho}
\eeq
and
\beq\label{8}
\lim_{R\to \infty} \HD(\esc_R (f)) \leq \frac{2M\rho}{2+M\rho}.
\eeq
\ethm

In this article we determine the exact values of these Hausdorff dimensions and we do this for the whole class of functions in Theorem \ref{BK theo}, and beyond.
The exact value of the Hausdorff dimension of the escaping set has been so far known only in some special cases, mainly for elliptic functions (see \cite{KUbook}), for the functions of the form $R(e^z)$, where $R$ is a rational function, obtained in the paper \cite{GaKo2016} by   Galazka and Kotus, and for Nevanlinna functions obtained in the paper \cite{Cui2019} by  Cui, where the special case $z\mapsto \l\tan z$ has been treated in \cite{GaKo2018}. 
Notice that for this later case one has strict inequality 
$$
\HD(\escf) < \frac{2M\rho}{2+M\rho}.
$$
Equality holds for all elliptic functions and for some interesting examples considered in  \cite{BwKo2012}.
Finally, there are also original new examples provided by Aspenberg and Cui in \cite{AspCui-I} which consist in some functions 
that have finite order and a finite singular set, that belong to Speiser class $\cS$, and whose escaping set may have any Hausdorff dimension in $[0,2]$.


Let $f$ be like in Theorem \ref{BK theo}, 
let $(a_j )_{j=1}^\infty$ be the sequence of all its poles, and let $m_j$ be the multiplicity of $a_j$, $j \ge 1$. Then 
\beq \label{2}
\limsup_{j\to\infty} m_j =M <\infty \quad \text{and}\quad  f(z)\sim\left(\frac{b_j}{z-a_j}\right)^{m_j} \quad \text{as}\quad z\to a_j
\eeq
for every $j\ge 1$ with some $b_j\in  \C\setminus \{0\}$.
Bergweiler and Kotus considered the series
\beq \label{3}
\sum_{j \; s.t.\;  a_j\neq 0} \left(\frac{|b_j|}{|a_j|^{1+1/M}} \right)^t,
\eeq
and an essential step in their proof of Theorem \ref{BK theo} was to show the convergence of it for all $t>\frac{2M\rho}{2+M\rho}$.

\bdfn\label{critical expo}
The critical exponent $\crex$ of the series  \eqref{3} is the infimum over all $t\in \R$ for which the
series \eqref{3} converges. 
\edfn

Notice that $\crex$ does depend on $M$ which itself does depend on $f$. Also, as explained above, this critical exponent  is at most $\frac{2M\rho}{2+M\rho}$
but in general it is strictly smaller. On the other hand, the functions under consideration have infinitely many poles and thus $\crex \geq 0$.

\smallskip

The main result of this article is the following result. Not only that it holds for all  functions of Theorem \ref{BK theo} but it also applies to functions of infinite order.

\bthm\label{theo main} 
Let $f$ be a meromorphic function in class $\cB$ for which $\infty$ is not an asymptotic value 
and whose poles have uniformly bounded multiplicities. Then
$$
\HD(\escf)= \lim_{R\to\infty}\HD(\esc_R(f))=  \crex,
$$
where $\crex$ is the critical exponent of the series \eqref{3}.
\ethm

Rempe-Gillen and Stallard called in \cite{ReSta2010} the number $\lim_{R\to\infty}\esc_R(f)$ the  \emph{eventual dimension} of $f$ and asked if it can differ from
 $\HD(\escf)$. Theorem \ref{theo main} shows that this is not the case for any function  considered therein.

Combining our result with the aforementioned results by Bara\'nski and Schubert, one obtains the following full description for the Hausdorff dimension of escaping sets.

\bcor
Let $f\in \cB$ be a meromorphic function of finite order with poles having bounded multiplicity. Then
\ben
\item $\HD(\escf) =2$ if  $\infty$ is an asymptotic value of $f$ and 

\,

\item $\HD(\escf)= \crex$ otherwise.
\een
\ecor

The upper bound in Theorem \ref{BK theo} does depend on the order $\rho$ of $f$ and converges to $2$ when $\rho \to \infty$.
This suggests that $\HD(\escf) =2$ whenever $f$ has infinite order, and this is the reason why in this case the behavior of $\escf$ with respect
to Lebesgue measure has been studied rather than the Hausdorff dimension of $\escf$ (see \cite[Theorem 1.3]{BwKo2012} and the discussion preceding it). However, our result, i.e. Theorem \ref{theo main}, shows that $\HD(\escf)$ does not really depend on the order. It certainly depends on the distribution of the poles but is also influenced 
by the coefficients $b_j$.

It is quite easy to determine the critical exponent $\crex$ as soon as the poles $a_j$ and the coefficients $b_j$ are known.
We exploit this in providing some quite simple explicit examples in Section \ref{s7}. Among them are functions having the following properties
which, up to our best knowledge, have not appeared in the literature.

\bthm\label{theo examples}
There exist meromorphic functions $f\in \cB$ of positive finite and also of infinite order for which $\infty$ is not an asymptotic value, whose all poles have bounded multiplicities, and  
$$
\HD(\escf)=0.
$$
\ethm

Our paper is organized as follows. First we present some background in Section \ref{s2}.
In the next section, Section \ref{s3}, we collect some properties of the behavior of the functions from Theorem~\ref{theo main}, near poles. 
The first step of the proof of Theorem \ref{theo main} is also done therein.
This part provides an upper estimate of $\HD(\escf)$ which is a
straightforward adaption of the proof from \cite{BwKo2012}. Indeed, 
Bergweiler and Kotus used the most natural dynamical covers of the escaping set, in fact of $\esc _R(f)$, 
from which the upper bound in Theorem \ref{BK theo} follows. The starting point of the present paper is the observation that their argument 
gives in fact that
$$
\hd(\escf)\leq  \lim_{R\to\infty} \HD(\esc _R(f))\leq \crex.
$$

Establishing a precise lower bound is the new and delicate step. For all specific functions for which $\HD(\escf)$ was known, such as elliptic and Nevanlinna functions mentioned above, the precise lower bound was obtained using  McMullen's method in \cite{McM87}. However, for this method to give sharp estimates the poles must be very 
uniformly distributed in addition with a uniform behavior of coefficients $b_j$, $j\ge 1$. This is clearly the case for elliptic and Nevanlinna functions as well
as for the examples in \cite{BwKo2012}.
In order to treat our general case we had to proceed totally differently. Our idea was to use the concept of non-autonomous iterated function systems and some results proved for them. Indeed, in Section\ref{s4}, we briefly summarize some facts from this theory. Then, in Section \ref{s5}, we prove the inequality 
$$
\hd(\escf)\geq  \crex,
$$ 
by associating to the escaping set $\escf$ some non-autonomous conformal iterated function system. This allows us to apply then some results of Rempe-Gillen and the second named author obtained in \cite{RU2016}. We construct our non-autonomous iterated function system
by carefully choosing appropriate inverse branches of the third iterate of the function $f$, and by doing some desingularizing change of variables which moves infinity to zero. Theorem~\ref{theo main} is then proved.

In Section~\ref{s6} we slightly generalize Theorem~\ref{theo main}. Our method being of local nature, we can consider arbitrary functions of class $\cB$
having infinitely many poles and determine the dimension of the points whose orbit escape through poles having multiplicity bounded by any given number $M$. We illustrate this by considering Euler Gamma functions which have in the same time a logarithmic tract and infinity many simple poles. 
From Bara\'nski's
\cite{B08} and Schubert's \cite{Schubert} work we know that the set of points of the escaping set whose orbit eventually are in the logarithmic tract has full dimension two. Here we get the additional information that the set of points escaping to infinity through poles is zero.

Finally in Section~\ref{s7} we prove Theorem~\ref{theo examples} by providing simple concrete constructions.

\section{Preliminaries} \label{s2}

The disk centered at $a\in \C$ and of radius $R>0$ will be denoted by $\D(a,R)$ and we define
$$
\D_R^*:=\{z\in \C \; : \; |z|>R\}.$$
If for some constant $C\geq 1$ we have $A/C\leq B\leq CA$ then this will be denoted by $A\asymp B$. If we need quantitative estimates, we write more precisely 
$$
A\asymp^C B.
$$
The inverse of a transcendental meromorphic function $f$ can have two types of singularities: asymptotic and critical values.
A point $b\in \cbar$ is an asymptotic value of $f$ if there  exists a curve $\ga \subset \C$ tending to infinity  such that $f(z)\to b$
as $z\to\infty$ along $\ga$. A critical value of $f$ is any point of the form $f(c)$ where 
$$
c\in \cri=\{z\in \C\,:\; f \text{ is defined at $z$ and } f'(z)\neq 0\},
$$
the set of critical points of $f$. The set of all finite singular values of $f$ is denoted by $\sing (f)$ and is commonly referred to as the singular set of $f$. A detailed description of singular values of meromorphic functions can be found for example in \cite{BwEr-Singularities}. 

We repeat from the introduction that class $\cB$ consists in all those meromorphic functions $f$ for which the singular set $\singf$ is bounded. This class goes back to
Eremenko and Lyubich \cite{EL92} who introduced it for entire functions. The functions in $\cB$ are also frequently called of \emph{bounded type}.

Unless stated otherwise, we consider in this paper transcendental meromorphic functions $f$ having the following properties:
\ben
\item[(i)] \label{item i} $f\in \cB$.
\item[(ii)] \label{item ii} $\infty$ is not asymptotic value of $f$.
\item[(iii)] \label{item iii} All poles $a_j$ of $f$ have uniformly bounded multiplicities denoted by $m_j$, so that the number $M$ in \eqref{2} is well defined, and consequently so is the critical exponent $\crex$. Also, we can and we will assume without loss of generality that $0$ is not a pole of $f$.
\een
 One should always have in mind that a transcendental meromorphic function $f$ for which $\infty$ is not an asymptotic value
 always  has infinitely many poles. Moreover, we can use them to characterize the order $\rho=\rho(f)$ of the function $f$ as the the critical exponent of the series
 \beq\label{38}
 \sum_{j=1}^\infty \frac{1}{|a_j|^t},
 \eeq
meaning that this series converges if $t>\rho$ and diverges if $t<\rho$; all these facts are in \cite[Section 2 ]{BwKo2012} and they go back to Iversen and Teichm\"uller.

 \smallskip

We fix a number $R_0>0$ such that 
\beq\label{11}
\sing (f) \subset \D (0,R_0/2) \  \text{ and } \  f(0)\subset \D (0,R_0/2).
\eeq
  Let $R\geq R_0$.
 Since $\infty$ is not an asymptotic value, Lemma 2.1 in \cite{RS99a} asserts that $f^{-1} (\D_R^* )$ is a countable union of simply connected bounded components 
 $U_j(R)$, $j \ge 1$, each of which contains exactly one pole $a_j\in U_j(R)$. When $R\geq R_0$ is fixed then we often denote these components simply  $U_j$.

\smallskip

All relevant information on the dynamics of transcendental functions, in particular the definition of the Julia set, can be found in the survey \cite{Bergweiler-survey} by Bergweiler, see also \cite{KUbook}. We deal here only with the escaping set $\escf$ and with the sets $\esc_R(f)$ which have been introduced in \eqref{5} and \eqref{6}.

Falconer's book \cite{Fal14}, \cite{PUbook}, and  \cite{KUbook} contain all relevant information on fractal dimensions we need; especially the two latter of them are designed for dynamical needs. One can find the definition of Hausdorff dimension in Section 3
and we  recall that the Hausdorff dimension of a set $E$ will be denoted by $\HD (E)$. Many results on fractal dimensions for meromorphic functions have been obtained by several
authors, the interested reader can find an overview of them in the survey  \cite{Stallard-Survey} by Stallard.

\smallskip

The restriction of every conformal map $g: \D (0, 2r) \to \C$ to $\D(0,r)$ has  \emph{bounded distortion} in the sense that there exists an absolute constant $K\geq 1$ such that 
\beq\label{bdd distortion}
|g'(z_1)|\leq K |g'(z_2)| \quad \text{for all $z_1,z_2\in \D(0,r)$.}
\eeq
This is  the well known Koebe's distortion theorem.

We will need  its more sophisticated version. It relies on the concept of modulus of ring domains that we briefly recall now,
see for example  \cite{Ahl06} for more details.

The modulus of the round annulus $A(r,R)=\{r<|z|<R\}$ is 
$$mod (A(r,R))= \frac 1{2\pi}\log \frac Rr \quad , \quad 0<r<R$$
and the modulus is a conformal invariant which can be used to define the modulus of general annuli.
These  are domains
$\cA\subset \C$ such that, for some radii $0<r_1<r_2$, $\cA$ is the conformal image of $A(r_1,r_2)$ and thus in fact the conformal image
of $A(1,R)$ where $R=\frac {r_2}{r_1}$. Set then 
$$
mod (\cA) =  \frac 1{2\pi}\log R \, .
$$
We need the following important fact which results from Theorem 2 in  \cite[p.10]{Ahl06}.
\bfact\label{fact 1}
If $\cA\subset \C$ is an annuli containing a round annuli
$A(r,R)\subset \cA$ for some $0<r<R<\infty$
then $$mod (\cA) \geq \mod A(r,R) =\frac 1{2\pi} \log R/r.$$
\efact
The following generalization of \eqref{bdd distortion} is contained in Theorem 2.9 in \cite{McMullen94}.

\bthm\label{koebe general}
For every $m>0$ there exists $K(m)\geq 1$ such that the following holds.
Let $U$ and $V$ be open domains in $\C$ both conformally equivalent to $\D$ such that $U\subset \ov U\subset V$ and 
$$
\mod(V\sms \ov U)\geq m.
$$
If $F:V\to \C$ is holomorphic univalent map. Then
$$
|F'(z_2)|\leq K(m) |F'(z_1)| \quad \text{for all $z_1,z_2\in U$.}
$$
\ethm

\section{Estimates and upper bound} \label{s3}
This part is a straightforward adaption 
of \cite{BwKo2012}. Throughout the whole section $f$ is a meromorphic function satisfying all the hypotheses of Theorem~\ref{theo main}.
We first give a fairly complete description of the behavior of the function $f$ near its poles and then we give the required upper estimate
of the Hausdorff dimension of $\esc_R(f)$ for large $R$.

\subsection{Estimates near poles} The behavior of a function $f$ near a pole $a_j$ follows from  \eqref{5}. So, let us recall it:
$$
f(z)\sim \left(\frac{b_j}{z-a_j}\right)^{m_j}
$$ 
as $z\to a_j$ for some complex number $b_j\neq 0$ and with 
$m_j$ being the multiplicity of $f$ at $a_j$ that satisfies  $1\leq m_j\leq M$ except for at most finitely many poles. Since we will only consider poles that lie outside some large disk we can assume that 
\beq\label{12}
1\leq m_i\leq M=\limsup_{j\to\infty} m_j.
\eeq
for every $i\ge 1$. We recall that $R_0>0$ has been fixed large enough  so that \eqref{11} holds. We shall prove the following.

\blem\label{14} 
For every $j\ge 1$ there exists a conformal homeomorphism $\Psi_j : U_j (R_0) \to \D(0, R_0^{-1/m_j} )$ such that $\Psi_j(a_j)=0$, $\Psi_j' (a_j)=1/b_j$, and 
$$
f|_{U_j (R_0)}= 1/ \Psi_j ^{m_j}
$$
\elem

\bpf 
Let $\Psi_j: U_j(R_0) \to \D(0, R_0^{-1/m_j} )$ be any conformal homeomorphism with $\Psi(a_j)=0$. It follows from the asymptotic behavior of $f$ near $a_j$ that the function $f\, \Psi_j^{m_j}$ is well defined, holomorphic, and non-zero on $U_j(R_0)$ and in addition  
$$
\lim_{U_j(R_0)\ni z\to \partial U_j(R_0)} | f(z)\,  \Psi_j^{m_j}(z) |=1.
$$
Thus, the maximum principle implies that 
$f\,  \Psi_j^{m_j}$ is constant and of modulus one. 
The required conformal map is now just  this $\Psi_j$ normalized by a rotation and then clearly  $\Psi_j' (a_j)=1/b_j$.
\epf

As the  first application of this lemma, we can now describe the geometry of the domains $U_j$. In the next lemma, the constant $K$ is the absolute constant from Koebe's distortion property \eqref{bdd distortion}.

\blem\label{15} 
For every $R\geq 2^MR_0$,
\beq \label{16}
\D\Big(a_j , \frac1K R^{-1/m_j}|b_j|\Big) \subset U_j(R) \subset \D\Big(a_j , K R^{-1/m_j}|b_j|\Big).
\eeq
Moreover, there exists $R_1\geq 2^MR_0$ such that for every $R\geq R_1$ there exists $\kappa_R > 1$ such that
$\kappa _{R}\le 2$ for every $R\ge R_1$, $\lim_{r\to\infty} \kappa _R =1$, and 
\beq\label{17}
\frac 1{\ka_R} \leq \frac{|z|}{|a_j|}\leq \ka_R
\eeq
for every $z\in U_j(R)$.
\elem

\bpf
We have $\Psi_j(U_j(R))=\D(0, R^{-1/m_j})$ where $\Psi_j $ is
 the conformal map of Lemma \ref{14}. Since $m_j\leq M$, we have that
 $R^{-1/m_j} \leq \frac12 R_0 ^{-1/m_j} $. This allows us to apply the bounded distortion property \eqref{bdd distortion} to $\Psi_j^{-1}$ from which \eqref{16} follows.

First let $R := 2^MR_0$. By \eqref{11}, $0\not\in U_j(R)$ and, using the first inclusion of \eqref{16}, we see that 
\beq\label{31}
|a_j| - \frac{|b_j| }{K R^{1/m_j}}>0, \quad \text{thus} \quad 
|b_j|\leq KR|a_j| = K2^MR_0 |a_j|.
\eeq
Now, take an arbitrary $R \geq 2^MR_0$ and then $z\in U_j(R)$. Then,
$$
|z-a_j|\leq \frac{K}{R^{1/m_j}}|b_j|\leq \frac{K}{R^{1/M}}  K2^MR_0 |a_j|,
$$
and thus
$$
\Big|\frac{z}{a_j}-1\Big|\leq 2^M K^2\frac{R_0}{R^{1/M}}$$
which implies \eqref{17} if $R$ is large enough, which determines $R_1$.
\epf

The next technical result provides  uniform estimates for the derivative $f'$ in the domains $U_j(R_1)$.

\blem\label{18} For every $j\ge 1$ we have that
$$
|f'| \asymp^{KM} \frac{1}{|b_j|}|f|^{1+\frac{1}{m_j}}  \quad \text{in $\;\;U_j(R_1)$.}
$$
\elem

\bpf
By Lemma \ref{14} we have $1/f=\Psi_j^{m_j}$ in $U_j$. Differentiating this relation yields
$$
|f'|=  m_j |\Psi_j^{m_j-1}\Psi_j' f^2|= m_j|\Psi_j'| |f|^{1+\frac{1}{m_j}}\asymp m_j \frac{1}{|b_j|}|f|^{1+\frac{1}{m_j}} \quad \text{in $\; U_j$}
$$
by \eqref{bdd distortion} and since $\Psi'(a_j)=1/b_j$. This proves the lemma since $m_j\leq M$.
\epf
 
\subsection{Upper estimate} Let us recall that this part is taken from \cite{BwKo2012} except that we alter the final conclusion.
We include it for the sake of completeness and as a preparation for the proof of the lower bound. Here is the precise statement.

\bprop\label{prop upper bd}
Let $f$ be a meromorphic function of class $\cB$ for which $\infty$ is not an asymptotic value 
and with poles having bounded multiplicity. Then
$$
\HD(\escf)\leq  \lim_{R\to\infty}\HD(\esc_R(f))\leq  \crex,
$$
where $\crex$ is the critical exponent of the series \eqref{3}. 
\eprop

\bpf
For every $R\geq R_1$ define the set
$$
\esc_R^0(f):=\big\{z\in\C:f^n(z)\ge R \quad \text{for all} \quad n\ge 0\big\}
$$
and, given an integer $l\ge 0$,
$$
\esc_{R,l}^0(f):=\big\{z\in \D^*_{R}\; : \; f^k(z)\in  \D^*_{R} \quad \text{for all} \quad 1\leq k<l\big\} \quad , \quad l\geq 0.
$$
Of course
$$
\esc_R^0(f)\subset \esc_{R,l}^0(f)
$$
for every $l\ge 0$. We first will estimate the Hausdorff dimension of the $\esc_R^0(f)$; for the ease of notation we will do this for the sets $\esc_{4R}^0(f)$

The pole neighborhoods $U_j$ and their successive pullbacks by inverse branches of $f$ give natural covers of the sets $\esc_{4R,l}^0(f)$, $l\ge 2$. We now explain it.
By Lemma \ref{15} we have $U_j\subset \D_R^*$ whenever
$U_j\cap \D_{4R}^*\neq \emptyset$ since then $|a_j| > 2 R$, whence 
$$
|z|\geq \frac 12 |a_j| > R \quad \text{for every} \quad z\in U_j.
$$
Since $\singf \cap \D^*_{R_0}=\emptyset$, all holomorphic inverse branches of $f$ are well defined on every simply connected
subdomain of $\D^*_{R_0}$. In particular, given any $j,k$ such that $U_j,U_k\subset \D^*_{R_0}$, there exist precisely $m_j$ such 
branches that map $U_k$ into $U_j$. Let
 $g_{j,k}:U_k \to U_j$ be any of them. Lemma \ref{18} along with \eqref{17} give precise estimate for its derivative.
  \beq\label{21}
 |g'_{j,k}(z) | \asymp^{KM} |b_j|/|z|^{1+\frac{1}{m_j}}\asymp^2  |b_j|/|a_k|^{1+\frac{1}{m_j}}\leq  |b_j|/|a_k|^{1+\frac{1}{M}}
 \;\; , \;\; z\in U_k.
 \eeq
Given $l\ge 2$ denote by $Z_l(R)$ the set of all $l$--tuples $\omega=(j_1,...,j_l)$ of positive integers such that $U_{j_1},..., U_{j_l}\subset \D_R^*$. Then set
\beq\label{20}
 V_\om:=g_{j_1,j_2}\circ ...\circ g_{j_{l-1},j_l} (U_{j_l}),
 \eeq
and we have 
$$
\bigcup_{\om\in Z_l}V_\om\spt \esc_{4R,l}^0(f)\spt \esc_{4R}^0(f).
$$
We have to estimate from above spherical diameters of the sets $V_\om$, $\om\in Z_l(R)$, in order to get an upper bound of $\HD\(\esc_{4R}^0(f)\)$. Notice that 
 $$
 \diam_{sph} (V)\asymp \diam (V)/ |a_{j_1}|^2
 \leq \diam (V)/ |a_{j_1}|^{1+\frac 1M}.
 $$
We are thus left to estimate the Euclidean diameter $\diam (V_\om)$. It follows from \eqref{16} that 
 $$
 \diam(U_{j_l})\leq K R^{-1/M}|b_{j_l}|.
 $$
 Employing \eqref{21} shows that there exists a constant $C$ such that
 $$
\diam(V_\om)\leq C^{l-1}\frac{|b_{j_1}|}{|a_{j_2}|^{1+\frac1M}}...\frac{|b_{j_{l-1}}|}{|a_{j_l}|^{1+\frac1M}}\diam(U_{j_l})
 $$
 and thus for the spherical derivative one gets
 $$
 \diam_{sph}(V_\om)\leq \frac {KC^{l-1}}{R^{1/M}} \prod_{k=1}^l \frac{|b_{j_k}|}{|a_{j_k}|^{1+\frac1M}}.
 $$
Summing now over all sets $V_\om$, $\om\in Z_l(R)$,  and taking into account that for each pair $j,k$ there are at most $M$ inverse branches of $f$ mapping $U_j$ into $U_k$, gives us for every $t>0$ the following estimate
 $$
 \sum_{\om\in Z_l(R)} \diam_{sph}(V)^t\leq \frac{K^t}{R^{t/M}} \left[ CM \sum_{j\in Z_0}
 \left(\frac{|b_{j}|}{|a_{j}|^{1+\frac1M}}\right)^t \right]^l
 $$
If $t>\crex$, the critical exponent of the series \eqref{3}, then for all $R>R_1$ large enough, 
$$
\sum_{\om\in Z_l(R)} \diam_{sph}(V_\om)^t\leq  \frac{K^t}{R^{t/M}} \frac1{2^l} \longrightarrow 0  \quad \text{as} \quad l\to \infty.
$$
Therefore,
$$
\HD\(\esc_{4R}^0(f)\)\le t.
$$
Hence 
$$
\HD\(f^{-k}\(\esc_{4R}^0(f)\)\)\le t
$$
for every integer $k\ge 0$. Since also
$$
\esc_{5R}(f)\sbt \bu_{k=0}^\infty f^{-k}\(\esc_{4R}^0(f)\),
$$
we thus get that
$$
\HD\(\esc_{5R}(f)\)\le t
$$
for all $R>R_1$ large enough. Thus, $\lim_{R\to\infty}\HD(\esc_R(f))\leq t$. Since $t>\crex$ was arbitrary, Proposition \ref{prop upper bd} follows. 
\epf
 
\brem\label{33} By \eqref{21}, in fact by Lemma \ref{18} along with \eqref{17}, we have for all $j, k\ge 1$ that
$$|g'_{j,k}(z) | \leq  L |b_j|/|a_k|^{1+\frac{1}{M}} \;\; , \;\;z\in U_k\;\; , \;\; \text{where}\;\;L=2MK.$$
On the other hand $|b_j|\leq K2^MR_0|a_j|$ (see \eqref{31}). Therefore,
$$
\text{$|g'_{j,k} | \leq  L K2^MR_0 \frac{1}{|a_k|^{1/M}}\leq  L K2^MR_0 \frac{1}{R_1^{1/M}}$ \ on $U_k$}.
$$

In the following we assume $R_1$ to be so large such that Lemma \ref{15} holds and that
\beq\label{32}
R_1^{1/M}\geq L K2^{M+1}R_0
\eeq
so that then 
$$
|g'_{j,k}(z)| \leq 1/2
$$
for all $j,k\ge 1$, meaning that all these maps are uniformly contracting.
\erem


\section{Non autonomous IFS} \label{s4}

A conformal iterated function system, IFS for short, in the complex plane is given by a family $\Phi=\{\Phi_i\}_{i\in I}$ of conformal contractions
$\Phi_i: \Om\to \Om$, where $I$ is an arbitrary countable set and where
$\Om$ is a bounded simply connected domain in $\C$ having some geometric properties that we will describe below. 
We assume that the functions $\Phi_i$ are in fact defined and continuous on $\ov \Om$ and that $\Phi_i (\ov \Om ) \subset \Om$.

One is interested by the dynamics of 
all possible compositions
$$
\Phi_{i_1}\circ ...\circ \Phi_{i_n}:\Om\longrightarrow\Om\;\; ,\quad i_k\in I \; \text{  for every } \; 1\leq k\leq n.
$$
Non-autonomous dynamical systems vary in time. In the setting of IFS this means that
for each $k$ the functions $\Phi_{i_k}$ in the above composition belong to a different family of functions. For example, consider the family $\Phi=\{\Phi_i\}_{i\in I}$
and define for each $n\geq 1$ a set 
$
I^{(n)}\subset I 
$ hence a family $\Phi ^{(n)} =\{\Phi_i\}_{i\in I^{(n)}}$.
The corresponding non autonomous dynamical system is given by the compositions
$$
\Phi_{i_1}\circ ...\circ \Phi_{i_n}:\Om\longrightarrow\Om\; ,\quad (i_1,...,i_n)\in \II_n=I^{(1)}\times ...\times I^{(n)} \; , \quad n\geq 1,
$$
and the limit set of this non-autonomous IFS is $\Phi ^{(n)} =\{\Phi_i\}_{i\in I^{(n)}}$, $n\geq 1$, is defined to be
$$
\J:= \bigcap_{n\geq 1}\bigcup_{(i_1,...,i_n)\in \II_n} \Phi_{i_1}\circ ...\circ \Phi_{i_n} (\Om).
$$
We only consider the case where all the sets $I^{(n)}$ are finite but their cardinalities are allowed to vary and grow with $n$, however only subexponentially, meaning that
\beq\label{40}
\lim_{n\to \infty } \frac 1n \log \# I^{(n)} =0 \,.
\eeq
We then say that the non--autonomous IFS $\Phi ^{(n)} =\{\Phi_i\}_{i\in I^{(n)}}$, $n\geq 1$, is of subexponential growth.
We are interested in the Hausdorff dimension $\HD(\J)$ of the limit set $\J$ of this system, and it is shown in \cite{RU2016} that this dimension is given by
a Bowen's Formula provided the system $\Phi$ has the  following properties.
\ben
\item $\Om$ is bounded and convex.
\,

\item Open Set Condition, meaning that 
$$
\Phi_i(\Om)\cap \Phi_j(\Om)=\es
$$
for all $n\ge 1$ and all $i, j\in I^{(n)}$.

\,

\item Bounded Distortion: all the compositions $\Phi_{i_1}\circ ...\circ \Phi_{i_n}:\Om\lra \Om$, $(i_1,...,i_n)\in \II_n$ have the property \eqref{bdd distortion} on $\Om$.

\, 

\item Uniform contraction: there exists $s\in (0,1)$ such that 
$$\text{$|\Phi '(z)|\leq s$ for every $z\in \Om$ and $i\in I$.}$$

\een

Also Bowen's Formula is not true in general but we have from Theorem 1.1 in \cite{RU2016}:

\bthm[Bowen's Formula]\label{bowen}
If a non--autonomous IFS $\Phi ^{(n)} =\{\Phi_i\}_{i\in I^{(n)}}$, $n\geq 1$, has all the properties listed above and is of subexponential growth, then 
$$
\HD(\J) = \inf \{t\geq 0 \; : \;\; \underline \P (t) <0\},
$$
where 
$$
\underline \P (t):= \liminf_{n\to\infty}\frac 1n \log \sum_{(i_1,...,i_n)\in \II_n}\| (\Phi_{i_1}\circ ...\circ \Phi_{i_n} )'\|_\infty ^t$$
is the, commonly called, lower pressure at $t$ of the non--autonomous system $\Phi ^{(n)} =\{\Phi_i\}_{i\in I^{(n)}}$, $n\geq 1$.
\ethm

\section{Lower bound} \label{s5}
Given Proposition \ref{prop upper bd}, Theorem \ref{theo main} will be proved if we can show that
$$
\HD(\escf)\geq \crex.
$$
This will be done here and the strategy is to associate to $f$ a non--autonomous 
IFS consisting of compositions of some selected holomorphic inverse branches of $f$,
for which Theorem \ref{bowen} (Bowen's Formula) will be applicable. First of all we have to determine a good domain $\Om$ on which our non--autonomous IFS is to act.

\subsection{The Domain for the IFS} \label{s51} For $0<\al_1 <\al_2< 2\pi$ let
$$
\cC (\al_1,\al_2): =\exp \lt(\big\{z\in \C \; : \; \al_1<\Im z<\al_2 \big\}\rt)
$$
be the cone delimited by the directions $\al_1,\al_2$. We recall that $\crex$ is the critical exponent of \eqref{3}.

\blem \label{28}
There exists a cone $\cC =\cC (\al , \al + 2\pi \, 3/4)$  such that 
$$\sum_{j \; s.t.\;  a_j\in \cC} \left(\frac{|b_j|}{|a_j|^{1+1/M}} \right)^t =\infty$$
for all $0\leq t < \crex$ and such that $\cC$ contains infinitely many poles $a_j$ of maximal multiplicity $m_j=M$.
\elem

\bpf
Since $\C\setminus \{0\}$ can be covered by five cones $\cC (\al_1,\al_2 ) $ with $(\al_2-\al_1)/2\pi= 1/4$, at least one of them
must contain infinitely many poles of maximal multiplicity. Denote it by $\cC_0$. 

The plane minus the origin can also be covered by two cones $\cC_1,\cC_2$ of the form 
$\cC (\al , \al + 2\pi \, 3/4)$ each of them containing $\cC_0$ so that both $\cC_1$ and $\cC_2$
contain infinitely many poles of maximal multiplicity. Now, let $t_0<\crex$. Then the series in Lemma \ref{27}
diverges if $\cC=\cC_1$ or if $\cC=\cC_2$. Furthermore, if this is the case, say for  $\cC=\cC_1$, then the
series with  $\cC=\cC_1$ diverges  for all $t\leq t_0$. The conclusion now comes by considering 
a sequence  $t_n\to \crex$ with all $t_n<\crex$.
\epf

For $R \geq  R_1$ let $$\cC_R:=\cC \cap \D_R^*$$ where $\cC=\cC (\al , \al + 2\pi \, 3/4)$ is the cone of Lemma \ref{28}.

\blem\label{29} 
There exists $R_2\geq  2K^{2M} R_1$ such that for every pole $a_j \in \cC_{2R_2}$ we have that
$$
U_j(R_2)\subset U_j\lt(\frac{R_2}{2K^{2M}}\rt) \subset  \hat\Om = \cC_{R_2}(\al- \pi /8, \al + 2\pi \,13/16).
$$
\elem

\bpf
This follows immediately from the inequality \eqref{17} in Lemma \ref{15}.
\epf

In the following, $R_2$ will be the number from this lemma and let
$$
\text{$I$ be the set of all $j$ such that the pole $a_j\in \cC_{2R_2}$.}
$$
Since $\hat \Om$ is a simply connected subdomain of $\D_{R_2}^*$ and $\D_{R_2}^*$ is disjoint from $S(f)$, all holomorphic inverse branches of $f$ are well defined on $\hat \Om$
and, for every $j\in I$, there exists $m_j$ such branches mapping $\hat \Om$ into $U_j= U_j(R_2)\subset \hat \Om$. Denote by
$$f_j^{-1}: \hat \Om \to U_j$$
one, arbitrarily chosen, of these inverse branches. We already gave a precise expression for the derivative of this function in the formula \eqref{21}. Here we repeat it. 
\beq\label{27}
|(f_j^{-1})'(z)|\asymp ^{KM} \frac{|b_j|}{|z|^{1+\frac{1}{m_j}}} \quad , \quad z\in \hat \Om.
\eeq

Since, by Lemma~\ref{29}, $f_j^{-1}(\hat \Om)\sbt \hat \Om$ for all $j\in I$,
we could build  an IFS out of this family of functions. However it would fail to have some necessary properties. Especially, given the above expression
for the derivative of $f_j^{-1}$, these functions have no bounded distortion property on $\hat \Om$. Nevertheless, there is a somewhat miraculous solution 
to this problem by applying, what we would call, a singular change of variables.

\subsection{Desingularization and Derivative Estimates} \label{s52} 

Consider the map $u:\cbar\to\cbar$ given by the formula
$$
u(\xi)= \xi ^{-M}.
$$
It maps  $\D\lt(0, R_2^{-\frac 1M}\rt)\setminus\{0\}$ onto $\D_{R_2}^*$ and admits all holomorphic inverse branches on $\hat \Om$. Choose again arbitrarily
one of them, denote it by $u_*^{-1}$, and define
$$\Om := u_*^{-1} (\hat \Om).$$
Notice that $\Om$ is bounded and convex.

Transfer then the maps $f_j^{-1}$ to $\Om$ by conjugation to get
$$
\ph_j := u_*^{-1}\circ f_j^{-1} \circ u : \Om \lra V_j := u_*^{-1} (U_j) \subset \Om \quad \text{for all}\quad j\in I.
$$
After having done this change of variables, not all new maps will have bounded distortion yet. However, the maps $\ph_j$ associated to poles with maximal multiplicity will
do have this property. The following lemma along with Lemma~\ref{30} are for us the two main ingredients to get, in Lemma~\ref{30'}, the required bounded distortion.

\blem\label{26} For every $j\in I$, we have that
$$
|\ph_j'(\xi)| 
\asymp^{2KM} \frac{|b_j|}{|a_j|^{1+\frac 1M}}|z|^{\frac{1}{M} -\frac 1{m_j}} 
=\frac{|b_j|}{|a_j|^{1+\frac 1M}}|\xi|^{\frac{M}{m_j} -1},  \;\xi\in \Om \text{ and } z=u(\xi).
$$
In particular, when $a_j$ is a pole with maximal multiplicity  $m_j=M$ then
$$
|\ph_j'(\xi)|\asymp^{2KM} \frac{|b_j|}{|a_j|^{1+\frac{1}{M}}} \ , \  \xi\in \Om .
$$
\elem

\bpf Let $\xi\in \Om$ and set $z:=u(\xi)=\xi^{-M}$. since $u\circ \ph_j=f_j^{-1}\circ u$, we get 
$$(\ph_j(\xi))^{-M}= f_j^{-1} (\xi^{-M}).$$
So, taking derivatives, gives
$$
|\ph_j'(\xi)|=|(f_j^{-1})'(z)|\frac{|\ph_j(\xi)|^{M+1}}{|\xi|^{M+1}}= |(f_j^{-1})'(z)|\frac{|z|^{1+\frac 1M}}{ |f_j^{-1}(z)|^{1+\frac 1M}}
$$
Since $f_j^{-1}(z) \in U_j$, Lemma \ref{15} implies that $|f_j^{-1}(z)|\asymp^2 |a_j|$. The affirmations in Lemma \ref{26} now result from \eqref{27}.
\epf

\subsection{Non-autonomous IFS} \label{s53}
Our goal now is to build an appropriate non--autonomous IFS out of the functions $\ph_j$.
 In order to do this, we first define an appropriate family of holomorphic maps on $\Om$ and verify that they have all necessary properties listed in Section~\ref{s4}. 
Then we explain how to choose for each $n\geq 1$ its subfamily, depending on $n$, so that all these subfamilies taken together define a conformal non--autonomous IFS that fulfills all the requirements of Bowen's Formula in Theorem \ref{bowen}.

There are several constraints for a construction of such IFS. First of all, it must consist in maps that are uniformly contracting and have bounded
distortion. Secondly, its limit set must be, up to the change of variables $u$,
a subset of the escaping set of $f$, and it must have a sufficiently large Hausdorff dimension. It turns out that the multiplicity of the poles play a crucial role and
we get limit sets of large dimension  only if the dynamics of the IFS pass sufficiently often through poles of maximal multiplicity. Since these poles are so important we introduce a second notation for them. Let $(j_l)_{l=1}^\infty$ be a strictly increasing sequence of elements in $I$ such that $(a_{j_l})_{l=1}^\infty$ are all the poles with multiplicity $m_{j_l}=M$ in the set $\{a_i\}_{i\in I}$. Define then
$$
\alm_{l,M}:=a_{j_l}, \, \blm_{l,M}:=b_{j_l}, \, \Vlm_{l,M}:=u^{-1}(\Ulm_{l,M})=u^{-1}(U_{j_l}) \text{ and } \phlm _{l,M}:=\ph_{j_l}.
$$
Also, since $I$ is infinite, we can assume without loss of generality that $I=\N$.

The maps $\ph_j$ cannot be used directly but it turns out that there is a way to
  compose these maps carefully in order to obtain a good family of generators for an IFS we want. 
  So, consider the following maps that correspond to inverse branches of the third iterate of $f$:
  \beq\label{ifs maps}
  \Phi_{l,j}:= \phlm _{l,M}\circ \ph_{j}\circ \phlm _{l,M}: \Om \lra \Om \quad, \quad l\geq 1 \; , \; j\in I.
  \eeq
  These have indeed bounded distortion on $\Om$ and are uniformly contracting:

\blem\label{30}
For every $j\in I$, every $l\ge 1$, and every $z\in\Om$, we have 
$$
 | \Phi'_{l,j}(z)| \asymp \frac{|\blm_{l,M}|^2}{|\alm_{l,M}|^{2+\frac{1}{m_j}+\frac{1}{M}}}\,
 \frac{|b_j|}{|a_j|^{1+\frac 1M}} 
 $$
  and  
  $$
   | \Phi'_{l,j}(z)|\leq \frac 14.
  $$
\elem

\bpf First of all
$$ 
\Phi'_{l,j}= \Big( \phlm _{l,M}\circ \ph_{j}\Big)'\circ  \phlm _{l,M} \;  \phlm ' _{l,M}.
$$ 
Lemma \ref{26} and the fact that  
$\alm  _{l,M}$ is a pole of maximal multiplicity yield
\beq\label{34}
|\phlm ' _{l,M}|\asymp^{2KM} \frac{|\blm_{l,M}|}{|\alm_{l,M}|^{1+\frac{1}{M}}} \quad \text{on}\quad \Om\,.
\eeq
Let $\xi\in \Om$ and $\eta=  \phlm _{l,M}(\xi)\in \Vlm_{l,M}$. We have to estimate the other factor of $ \Phi'_{l,j}$
at $\eta$. Since $u(\eta)\in \Ulm_{l,M}$, we have that
$$
|\eta|^{-M} =|u(\eta)|\asymp^2 |\alm_{l,M}|.
$$
Using now Lemma \ref{26} once more and remembering that $\alm_{l,M}$ has maximal multiplicity, we get 
$$
\begin{array}{rl}
 \Big|\Big( \phlm _{l,M}\circ \ph_{j}\Big)'(\eta)\Big|&\asymp^{(2KM)^2} \frac{|\blm_{l,M}|}{|\alm_{l,M}|^{1+\frac{1}{M}}}
 \frac{|b_j|}{|a_j|^{1+\frac{1}{M}}}|\eta|^{\frac{M}{m_j}-1}\\
& \asymp ^2 \frac{|\blm_{l,M}|}{|\alm_{l,M}|^{1+\frac{1}{M}}}\frac{|b_j|}{|a_j|^{1+\frac{1}{M}}}|\alm_{l,M}|^{\frac{1}{M}-\frac 1{m_j}}
= \frac{|\blm_{l,M}|}{|\alm_{l,M}|^{1+\frac{1}{m_j}}}\frac{|b_j|}{|a_j|^{1+\frac{1}{M}}}.
\end{array}
$$
Multiplying this with \eqref{34} gives the required estimate for $ |\Phi'_{l,j}|$.

It remains to verify the contraction property. The, already proved, first formula of our lemma gives
$$
  | \Phi'_{l,j}|\leq 2 L^3  \frac{|\blm_{l,M}|^2}{|\alm_{l,M}|^{2+\frac{1}{m_j}+\frac{1}{M}}}\,
 \frac{|b_j|}{|a_j|^{1+\frac 1M}}
$$
where $L=2KM$. By \eqref{31} we have $|b_i|\leq K2^MR_0|a_i|$ for all $i$. So, since $|a_i|> R_1$ and $m_i\leq M$, 
$$
| \Phi'_{l,j}|\leq 2 L^3 (K2^MR_0)^3 \frac{1}{R_1^{3/M}}.$$
Along with \eqref{32}, this implies that $ | \Phi'_{l,j}|\leq 1/4$. The proof is complete.
\epf

\smallskip

We also need bounded distortion for arbitrary compositions of the maps $\Phi_{l,j}$.
\blem\label{30'} There exists $\hat K\geq 1$ such that for all  $n\geq 2$ and all 
$(l_k,i_k)\in I$, $ 1\leq k\leq n$, the map
$$
F=\Phi_{l_1,i_1}\circ \Phi_{l_2,i_2}\circ ...\circ \Phi_{l_n,i_n}:\Om\lra\Om
$$
satisfies
\beq\label{30' dist}
|F'(z_2)|\leq \hat K |F '(z_1)| \quad \text{for all $\;\;z_1,z_2\in \Om$.}
\eeq
\elem

\smallskip

\bpf
Write $F=F_{n-1}\circ \Phi_{l_n,i_n} $ where $F_{n-1}=\Phi_{l_1,i_1}\circ ...\circ \Phi_{l_{n-1},i_{n-1}}$. 
The map $\Phi_{l_n,i_n}$ satisfies \eqref{30' dist} with some constant $\hat K$, say $\hat K_1\ge 1$, by Lemma \ref{30}. 
Therefore, Lemma \ref{30'} will hold if we can show that there exists $\hat K_2$ such that 
$$
|F'_{n-1}(z_1)|\leq \hat K_2 |F '_{n-1}(z_2)| 
$$
for all $n\geq 2$, all $(l_k,i_k)\in I$ with $ 1\leq k\leq n$, and all $z_1,z_2\in \Phi_{l_n,i_n} (\Om )$.

But this follows from the general Koebe distortion theorem, i.e. Theorem \ref{koebe general}, if we can find an open disk $D$ such that $\Phi_{l_n,i_n} (\Om )\subset D$, the map $F_{n-1}$ has a holomorphic univalent extension to $D$, and 
$$
\mod\(D \setminus \ov {\Phi_{l_n,i_n} (\Om ) }\)\geq m
$$
where $m>0$ is a constant independent of all $n\geq 1$, $(l_k,i_k)\in I$, and $ 1\leq k\leq n$.

We have that $ \Phi_{l_n,i_n} (\Om )\subset \Vlm_{l_n,M}=u_*^{-1}(\Ulm_{l_n,M})=u_*^{-1}(U_{j_{l_n}})$ and, by Lemma \ref{29},
 $$
 U_{j_{l_n}}(R_2)\subset  U_{j_{l_n}}(R_2/Q)\subset \hat \Om \quad \text{where}\quad Q= 2K^{2M}.
 $$
 The map $F_{n-1}$ is well defined, holomorphic, and univalent on the disk $$u_*^{-1}(U_{j_{l_n}}(R_2/Q)).$$ It will be our disk $D$, and as
$$
\mod\lt(u_*^{-1}(U_{j_{l_n}}(R_2/Q)\sms \ov{\Vlm_{l_n,M}}\rt)=\mod\lt( U_{j_{l_n}}(R_2/Q) \setminus \overline {U_{j_{l_n}}(R_2)}\rt),
$$
it remains to estimate from below the modulus of the annulus 
$$
\cA:=U_{j_{l_n}}(R_2/Q) \setminus \overline {U_{j_{l_n}}(R_2)}.
$$
This can be done by applying Lemma \ref{15}. Indeed, it follows from this lemma that 
 $$
 U_{j_{l_n}}(R_2)\subset  \D\Big(a_{j_{l_n}},KR_2^{-1/m_{j_{l_n}}} |b_{j_{l_n}}|\Big)$$
 and
 $$
 \D\Big(a_{j_{l_n}}, \frac 1K (R_2/Q)^{-1/m_{j_{l_n}}}  |b_{j_{l_n}}|\Big)\subset  U_{j_{l_n}}(R_2/Q).
 $$
 Thus, $\cA$ contains the round annulus
 $$
 A\Big( KR_2^{-1/m_{j_{l_n}}} |b_{j_{l_n}}|, \frac 1K (R_2/Q)^{-1/m_{j_{l_n}}}  |b_{j_{l_n}}| \Big),
 $$
 which implies, by Fact \ref{fact 1}, that 
 $$
 \mod (\cA)\geq \frac 1 {2\pi} \log \frac{Q^{1/m_{j_{l_n}}}}{K^2}\geq  \frac 1 {2\pi M}\log 2.
 $$
So, the proof is complete by putting $m:= \frac 1 {2\pi M}\log 2>0$ .
\epf

Having now the right family of contractions defined on an appropriate domain $\Om$  we are ready to define the
non-autonomous system.
In order to do so, we will choose for each $n\geq 1$ a set $I^{(n)}\subset \N\times I$ and then
consider compositions
$$
\Phi_{l_1,i_1}\circ \Phi_{l_2,i_2}\circ ...\circ \Phi_{l_n,i_n}:\Om\lra\Om
$$
with $(l_k,i_k)\subset I^{(k)}$ for all $1\leq k\leq n$. The choice of these index sets does depend on the parameter $t$. So, fix $t< \crex$ arbitrary. Then, by the definition of $\crex$, the series \eqref{3} is divergent for this $t$. Lemma \ref{30} then implies that 
$$
\sum_{j\geq j_l} \inf|\Phi'_{l,j}|^t =+\infty 
$$ 
for every $l\geq 1$, where $(j_l)_{l=1}^\infty$ is again the sequence corresponding to the poles of maximal multiplicity
$\alm_{l,M}= a_{j_l}$. Therefore, for every $l\geq 1$, there exists  $k_l> j_l$ with 
$$
 \sum_{j\geq j_l}^{k_l} \inf|\Phi'_{l,j}|^t \geq 2 
$$
and, adding some terms in these sums if necessary, we may assume that 
$$
k_{l+1}-j_{l+1}\geq k_l-j_l 
$$
for all $l\geq 1$.
Define now
$$
\hat I^{(l)}:= \{(l,j)\; : \; j_l\leq j\leq k_l\}, \  \  l\ge 1.
$$
$\# \hat I^{(l)}$, the number of elements of $I^{(l)}$, increases with $l$ but may vary a lot. This is the reason why we define the sets $I^{(n)}$ as follows. They interpolate 
the sets $\hat I^{(l)}$ but in such a way that 
$$
\# I^{(n+1)}\leq \# I^{(n)} +1
$$
for all $n\geq 1$. One way of doing this is the following. Let 
$$
I^{(1)}:= \hat I^{(1)}.
$$
If $\# \hat I^{(2)} = \# I^{(1)}$, then put 
$$
I^{(2)}:=\hat I^{(2)}.
$$
Otherwise, we have $\# \hat I^{(2)} = \# I^{(1)}+N_1$ for some $N_1\ge 1$. We then define inductively
$$
I^{(2)}:= I^{(1)}\cup\{(1,k_1+1)\},\, I^{(3)}:= I^{(2)}\cup\{(1,k_1+2)\},....,
$$
$$
I^{(N_1+1)}:= I^{(N_1)}\cup\{(1,k_1+N_1)\},
$$
and then
$$
I^{(N_1+2)}:=\hat I^{(2)}.
$$
By continuing inductively this definition, we obtain a sequence of sets  $\(I^{(n)}\)_{n=1}^\infty$ having the following properties.
\blem\label{36} For every $t< \crex$ there exist $\(I^{(n)}\)_{n=1}^\infty$ a sequence of sets,  and $(l_n)_{n=1}^\infty$, a diverging sequence of positive integers, such that for every $n\geq 1$ the following hold:
\ben
\item 
$\displaystyle I^{(n)} =\{(l_n,j): j_{l_n}\leq j\leq h_n\} \quad \text{with  some}\quad h_n\geq k_{l_n}.$

\,

\item  $\displaystyle \sum_{(l_n,j)\in I^{(n)}} \inf| \Phi'_{l_n,j}|^t \geq 2 \;$.
\medskip
\item $\# I^{(n+1)}\leq \# I^{(n)} +1$.
\een
\elem
The last property immediately implies that the system has subexponential growth, i.e. it satisfies \eqref{40}, meaning that
\beq\label{120210218}
\lim_{n\to \infty } \frac 1n \log \# I^{(n)} =0 \,.
\eeq
whereas the second property yields positivity of the lower pressure for the parameter $t$:
\beq \label{41}
\underline \P (t)= \liminf_{n\to\infty}\frac 1n \log \sum_{(i_1,...,i_n)\in \II_n}\| (\Phi_{(l_1,i_1)}\circ ...\circ \Phi_{(l_n,i_n)} )'\|_\infty ^t
\ge \log 2>0.
\eeq

\subsection{Conclusion of the Proof of Theorem \ref{theo main}}

Theorem \ref{theo main} results from Proposition \ref{prop upper bd} and the next result.

\bprop\label{prop lower bd}
If $\J$ is the limit set of the non-autonomous conformal IFS $\Phi ^{(n)} =\{\Phi_i\}_{i\in I^{(n)}}$, $n\geq 1$, constructed in the previous subsection, then 
$$
u(\J) \subset \escf
$$ 
and
$$
\HD(\J)\geq \d.
$$
\eprop

\bpf
The first assertion of this proposition follows immediately from the construction of the non-autonomous IFS $\Phi ^{(n)} =\{\Phi_i\}_{i\in I^{(n)}}$ done in the previous subsection. 

The second assertion follows immediately from Lemma~\ref{30}, Lemma~\ref{30'}, Formula \eqref{120210218}, Formula \eqref{41}, and Bowen's Formula, i.e. Theorem~\ref{bowen}. 
\epf

\section{General functions of class $\cB$}\label{s6}

As mentioned in the Introduction, McMullen and various other authors who extended his work \cite{McM87} estimated the Hausdorff dimension
of the Julia set by considering the escaping set and, more precisely, the set of points whose orbit escapes through a tract of logarithmic singularity.
Here we consider, for general functions of class $\cB$, orbits that escape to infinity ``through poles''.

To make this more precise, let $f\in \cB$ be an arbitrary transcendental meromorphic function and let again $R_0>0$ be such that \eqref{11} holds
so that in particular $\singf \subset \D(0, R_0/2)$. For $R\geq R_0$ let again $U_j (R)$ be the connected component of $f^{-1}(\D_R^*)$ that contains a pole $a_j$. Consider then the set
$$
 \quad \cU_{M,R} :=\bigcup_{j \, :\, m_j \leq M}  U_j (R)
$$
where, as always, $m_j$ is the multiplicity of $a_j$ and where $M\geq 1$. Consider now the subsets 
$$
\esc (f,\cU_{M,R}) 
$$
of the escaping set of $f$ consisting of all those points $z\in \escf$ for which there exists $N\geq 1$ such that 
$$
f^n(z)\in \cU_{M,R}
$$ 
for all $n\geq N$. Clearly, these subsets of the escaping set are only 
relevant if $f$ has infinitely many poles. If $f$ is entire or has only finitely many poles  then the sets  $\esc (f,\cU_{M,R}) $ are empty.

\bthm\label{thm main extension}
Let $f\in\cB$ be a transcendental meromorphic function and assume that there exists $M_f$ such that  
$f$ has infinitely many poles $a_j$, $j\ge 1$, with multiplicities $m_j\leq M_f$. Then, for every $M\geq M_f$,
there exists $R_M >R_0$ such that for all $R\geq R_M$
$$
\HD\lt(\esc (f,\cU_{M,R})\rt) =\crex _M,
$$
where $\crex _M$ is the critical exponent of the series
$$
\sum_{j \; s.t.\;  m_j\leq M} \left(\frac{|b_j|}{|a_j|^{1+1/{\mu}}} \right)^t
$$
and where $\mu$ is the upper limit of all the numbers $m_j$ for which $m_j \leq M$.
\ethm

\bpf 
The proofs of the Propositions \ref{prop upper bd} and \ref{prop lower bd} are both 
local. Indeed, let $M\geq M_f$ and consider the restriction of $f$ to the set $\cU_{M,R_0}$.
All   inverse branches $f_j^{-1}:\D_{R_0}^* \to U_j (R_0)$, $m_j\leq M$, are well defined
and we can proceed exactly as in these proofs using this family of inverse branches in place of the full family of inverse branches considered in the proofs of the Propositions \ref{prop upper bd} and \ref{prop lower bd}.
\epf

Theorem \ref{thm main extension} contains  Theorem \ref{theo main} as its special case. Also, it is clear that it gives the following estimate valid for all functions in class $\cB$ having infinitely many poles. 
The critical exponents $\crex_M$ are increasing with $M$ so that $\sup_M \crex_M=\lim_{M\to \infty} \crex_M$. Denote this limit by $\crex_\infty$.

\bcor\label{c120210326}
If $f\in\cB$ is a transcendental meromorphic function having infinitely many poles, then
$$
\HD(\escf) \geq \crex_\infty \,.
$$
\ecor

Let us recall that the results of Bara\'nski (\cite{B08}) and Schubert (\cite{Schubert}) show that
$\HD(\escf )=2$ if $f\in \cB$ and if $\infty$ is asymptotic value. Therefore, Corollary \ref{c120210326}
 is meaningful for the functions of $\cB$ for which $\infty$ is not asymptotic value.

\smallskip

We finally provide a class of examples that illustrate 
Theorem \ref{thm main extension}. Let
$
\Gamma:\C\to\oc
$
be the Euler Gamma function and, for every $a\in\C$, let $\Gamma_a:\C\to\oc$ be the  translated  version of it defined by
$$
\Gamma_a(z)=\Gamma(z+a) \quad \text{, $\;z\in\C$. }
$$
It is well known and easy to check that $\Gamma$ and thus all the functions $\Gamma_a$, $a\in \C$,  are in class $\cB$
and have precisely one logarithmic tract over infinity. Moreover, the poles of a function $\Gamma_a$ are 
 $a_j=-j-a$, all these poles are simple and the residue of $\Gamma_a$ at $a_j$ is
 $$\text{$b_j =(-1)^j/j!\;\; $, $\;\;j\geq 1$. }$$

Since $\Gamma _a$ has a logarithmic tract over infinity, again \cite{B08} and \cite{Schubert} imply that 
\beq\label{220210326}
\HD(\esc (\Gamma _a)) =2.
\eeq
In fact, it follows from these papers that the Hausdorff dimension of the 
subset of the escaping set consisting of all points $z\in \esc (\Gamma _a)$
 whose orbits are eventually
contained in the tract  is equal to $2$.

Theorem \ref{thm main extension} allows to determine the dimension of the points that escape to infinity through the poles.
We recall that the later are simple poles so that $M_{\Gamma_a}=1$.

\blem
For all $R$ large enough $\HD \big(\esc (\Gamma _a ,  \cU_{1,R}) =0$.
\elem

\bpf
Let $R>0$ be such that Theorem \ref{thm main extension} applies. Then $\HD \big(\esc (\Gamma _a ,  \cU_{1,R})$ equals the critical 
exponent $\d_1$ of the series
$$
\sum_{j=1}^\infty \left(\frac{|b_j|}{|a_j|^{1+1/1}}\right)^t
\asymp \sum_{j=1}^\infty \left(\frac{1}{j!j^2}\right)^t
$$
which clearly is $\d_1=0$.
\epf

\section{Examples}\label{s7}

Using our main result, i.e. Theorem \ref{theo main}, it is particularly easy to determine the value of the Hausdorff dimension of the escaping set
as soon as one has a sufficiently good information about the poles of the function. We illustrate this here by providing some simple examples. The first class of examples consists in functions with 
finite order, in fact order two (but this can be easily generalized to any finite order), such that $\HD(\escf)$ realizes any value in $[0,2)$. Then we show, by explicit construction, that
functions of infinite order fulfilling the hypotheses of Theorem \ref{theo main} can have $\HD(\escf) <2$ and even $\HD(\escf)=0$.

Theorem \ref{theo examples} follows from the Theorems~\ref{prop ex1} and \ref{prop ex2} proven below. 

\subsection{Examples of Finite Order} 
We consider functions of the form
$$
f(z)= \sum_{|j|, |k| \geq N} \left( \frac{b_{j,k}}{z-a_{j,k}}\right)^M
$$
with poles $a_{j,k}= j+ik$, $j,k\in \Z$, all of same multiplicity $M\geq 1$. Remember that the neighborhood $U_{j,k}$ of a pole $a_{j,k}$ is approximately a disk
of size $|b_{j,k}|$; see Lemma \ref{15}. Hence, these numbers can not be too large and we take 
$$
0<b_{j,k} \leq |a_{j,k}|^{-\al} \quad \text{for some fixed $\al >2/M$}.
$$
Since $\al M>2$, we can fix the integer $N\geq 1$ sufficiently large so that 
\beq \label{37}
\sum_{|j|, |k| \geq N} |a_{j,k}|^{-\al M} \leq \frac 1{2^{M+1}}.
\eeq
Denote
$$
\Z_N^2:=\{(j,k)\in\Z^2:|j|, |k| \geq N\}
$$
\bthm\label{prop ex1}
Under these assumptions, $f$ is a well defined meromorphic function in class $\cB$ and of order two for which $\infty$ is not an asymptotic value and all of its poles have multiplicity $M$. furthermore:
\ben
\item If $b_{j,k} =|a_{j,k}|^{-\al}$ for all $(j,k)\in\Z_N^2$, then 
$$
\displaystyle \HD(\escf) = \frac{2}{\al + 1+\frac 1M}\in \lt(0, \frac{2}{ 1+\frac 3M}\rt).
$$
\item If $b_{j,k} =\exp (-|a_{j,k}|)$ for all $(j,k)\in\Z_N^2$, then $\HD(\escf) =0$.
\een
\ethm

\bpf
Fix $z\in \C$ distinct from all the poles. For all those $(j,k)\in \Z_N^2$ for which $|a_{j,k}|\geq 2|z|$ (keep in mind that the set of $j,k$ failing to satisfy this inequality is finite), we have
$$
 \left| \frac{b_{j,k}}{z-a_{j,k}}\right|^M\leq  \left| \frac{|a_{j,k}|^{-\al}}{|a_{j,k}|-|z|}\right|^M \leq 2^M |a_{j,k}|^{-(1+\al)M}.
$$
So, since $\al M\geq 2$, the series defining $f$ converges absolutely and uniformly on some neighborhood of each such point $z$ and, moreover, $f:\C\to \cbar$ is a well defined meromorphic function with poles being equal to the numbers $a_{j,k}$, $(j,k)\in \Z_N^2$. Furthermore, $f$
is of order two since the series \eqref{38} converges if and only if $t>2$. 

\sp Consider the disjoint disks $\D(a_{j,k}, b_{j,k})$, $(j,k)\in \Z_N^2$, and let 
$$
\cU:=\displaystyle  \bigcap_{(j,k)\in \Z_N^2} \D^c(a_{j,k}, b_{j,k}).
$$ 
\bclaim\label{c220210218}
The set $\cri$ of critical points of $f$ satisfies
$\cri \subset\cU$.
\eclaim

\bpf
Indeed, since
$$
f'(z)= -M   \sum_{(j,k)\in \Z_N^2} \frac{b_{j,k}^M}{(z-a_{j,k})^{M+1}}
$$
we have for all  $(j,k)\in \Z_N^2$, and all $\displaystyle z\in  \D(a_{j,k}, b_{j,k})\setminus\{a_{j,k}\}$ that
$$
\begin{array}{rl}
|f'(z)|&\geq M\left( \frac{b_{j,k}^M}{|z-a_{j,k}|^{M+1}}   -  \sum_{(u,v)\in\Z_N^2\sms\{(j,k)\}}  \frac{b_{u,v}^M}{|z-a_{u,v}|^{M+1}}  \right)\\
&\geq M\left(b_{j,k}^{-1}-  \sum_{(u,v)\in\Z_N^2\sms\{(j,k)\}} 2^{M+1}b_{u,v}^M \right)
\end{array}
$$
since $|z-a_{u,v}|\geq 1/2$ if $(u,v)\neq (j,k)$. Thus,
$$
|f'(z)|\geq M\Big( |a_{j,k}|^\al  -  2^{M+1}\sum_{(u,v)\neq (j,k)}  |a_{u,v}|^{-\al M}  \Big) \geq
M\left( (\sqrt{2}N)^\al -1\right)>0.
$$
\epf
\smallskip

\bclaim\label{c120210218}
$|f(z)|\leq 2$ for every $z\in\cU$.
\eclaim

\bpf
Let $z\in \cU$. If there exists $(j,k)\in\Z_N^2$ such that $|z-a_{j,k}|\leq 1/2$, then $|z-a_{u,v}|\geq 1/2$ for all $(u,v)\in\Z_N^2\sms\{(j,k)\}$ and 
$$
\begin{array}{rl}
|f(z)|&\leq  \left(\frac{b_{j,k}}{|z-a_{j,k}|}\right)^M + \sum_{(u,v)\in\Z_N^2\sms\{(j,k)\}} \left(\frac{b_{u,v}}{|z-a_{u,v}|}\right)^M \\
&\leq  \left(\frac{b_{j,k}}{b_{j,k}}\right)^M + \sum_{(u,v)\in\Z_N^2\sms\{(j,k)\}} \left(\frac{b_{u,v}}{1/2}\right)^M \\
&\leq 1+ 2^M 2^{-(M+1)} \leq 2
\end{array}
$$
since $b_{u,v}\leq |a_{u,v}|^{-\al}$ and since we have \eqref{37}.

In the other case, i.e. if $|z-a_{j,k}| >1/2$ for all $(j,k)\in\Z_N^2$, then
$$
|f(z)| \leq  \sum_{(j,k)\in\Z_N^2} \left( \frac{b_{j,k}}{1/2}\right)^M \leq 2^M 2^{-(M+1)} \leq 2
$$
also since $b_{j,k}\leq |a_{j,k}|^{-\al}$ and since we have \eqref{37}. Thus the claim holds.
\epf

\bclaim
$$
\sing (f) \subset \ov \D(0,2),
$$
and consequently $f\in \cB$ and $\infty$ is not asymptotic value.
\eclaim

\bpf
It directly follows from Claims~\ref{c220210218} and \ref{c120210218} that
$$
f(\cri)\subset \ov \D (0,2).
$$
Concerning the asymptotic values, if $\ga \subset \C$ is a curve tending to infinity whose image under $f$ converges to an asymptotic value $b\in \hat \C$, then there are arbitrary large values $z\in \ga $ that are in $\cU$. Claim~\ref{c120210218} implies that for such values $z$ we have $|f(z)|\leq 2$, consequently $|b|\leq 2$. All of this shows that $\sing (f) \subset \ov \D(0,2)$ and that  $\infty$ is not an asymptotic value. We are done.
\epf

Thus the function $f$ satisfies all hypotheses of our Theorem \ref{theo main}, and thus it suffices to determine the critical exponent $\crex$ 
of \eqref{3} in order to  determine $\HD(\escf)$.

\

\noindent
\emph{Case (1):} $b_{j,k} =|a_{j,k}|^{-\al}$ for all $(j,k)\in\Z_N^2$. Then the series
$$
 \sum_{(j,k)\in\Z_N^2} \left( \frac{b_{j,k}}{|a_{j,k}|^{1+\frac{1}{M}}}\right)^t =
  \sum_{(j,k)\in\Z_N^2} |a_{j,k}|^{-\lt(\al +1+\frac{1}{M}\rt)t}
$$
converges if and only if $t\lt(\al +1+\frac{1}{M}\rt)> 2$. This, with the help of Theorem \ref{theo main}, shows that 
$$
\HD(\escf) =\crex = \frac{2}{\al +1+\frac{1}{M}}.
$$

\

\noindent
\emph{Case (2):} $b_{j,k} =\exp(-|a_{j,k}|)$  for all $(j,k)\in\Z_N^2$. Then the series
$$
 \sum_{(j,k)\in\Z_N^2} \left( \frac{b_{j,k}}{|a_{j,k}|^{1+\frac{1}{M}}}\right)^t =
 \sum_{(j,k)\in\Z_N^2} \left( \frac{\exp(-|a_{j,k}|)}{|a_{j,k}|^{1+\frac{1}{M}}}\right)^t
$$
converges for all $t>0$ which, in the same way as in Case~1 shows that 
$$
\HD(\escf) =\crex =0.
$$
\epf

\brem In both cases, i.e. if  $b_{j,k} =|a_{j,k}|^{-\al}$ or if $b_{j,k} =\exp(-|a_{j,k}|)$, we see that the series \eqref{3} diverges at the critical exponent $t=\crex$.
This in fact gives further information on the dimension $\dim \esc_R (f)$ since it was observed in \cite{BwKo2012} that one can then employ 
\cite{My09} in order to get
$$\HD(\esc_R (f)) > \crex $$
for every $R>0$ although, by Theorem \ref{theo main}, $\lim_{R\to\infty}\dim \esc_R (f) = \crex $.
\erem

\smallskip

\subsection{Examples of Infinite Order}
There are many ways to modify the preceding example in order to get such functions having infinite order. Here is one way to do it. Let
$$
a_j := \log j \, \text{,}\quad b_j :=e^{-j},
$$
and define 
\beq\label{120210222}
f(z):=\sum_{j\geq 8} \frac{b_j}{z-a_j}, \  \ z\in\C.
\eeq

\bthm\label{prop ex2}
The function $f$ defined above in the formula \eqref{120210222} is a well defined meromorphic function of class $\cB$ of infinite order infinity for which $\infty$ is not an asymptotic value.
Furthermore, all poles of $f$ have multiplicity one and  the escaping set of $f$ satisfies
$$
\HD(\escf) =0.
$$
\ethm

\bpf Since $\displaystyle \frac{b_j}{a_j}= \frac{e^{-j}}{\log j}$, $j\geq 8$, the same argument, even easier, as in the proof of Theorem~\ref{prop ex1}, shows that the series of \eqref{120210222} defines a meromophic function from $\C$ to $\cbar$ with poles being equal to the numbers $a_{j}$, $j\ge 8$. Furthermore, since $a_j = \log j$, $j\ge 8$, the series  \eqref{38} associated with our present function $f$ diverges for all $t>0$. This shows that the order $\rho$ of $f$ is equal to $\infty$.

The other estimates needed to prove Theorem~\ref{prop ex2} are (also) similar and even simpler than in the proof of Theorem~\ref{prop ex1}. Consider again the sets
$$
\D(a_{j},b_j), \ j\ge 8, \  \  {\rm and}  \  \  \cU:=\bigcap_{i \geq 8} \D^c(a_{i},b_i).
$$

Let $z\in \cU$. If $z\in \D\lt(a_j,\frac1{3j}\rt)$ for some $j\geq 8$ then, by an elementary calculation, 
\beq \label{39}
|z-a_k|\geq \frac1{3k} \  \  {\rm for \ all}\  \  k\neq j, k\geq 8.
\eeq
This implies 
$$
|f(z)| \leq  \frac{b_j}{|z-a_j|} + \sum_{k\neq j,\, k\geq 8}\frac{e^{-k}}{|z-a_k|} \leq 1+ 3\sum_{ k\geq 8} k e^{-k} \leq 2.
$$
If $z\not\in \D\lt(a_j,\frac1{3j}\rt)$ for all $j\geq 8$, then
$$
|f(z)| \leq   \sum_{j\geq 8}\frac{b_j}{|z-a_j|} \leq 3\sum_{ j\geq 8} je^{-j} \leq 1.
$$
This shows that 
$|f(z)|\leq 2$ for every $ z\in \cU$.

The next step is to show that $\cri \subset \cU$. We have 
$$
f'(z)= -\sum_{j\geq 8} \frac{b_j}{(z-a_j)^2}.
$$ 
Let $z\in\cU$. Then $z\in \D(a_{j},b_j)\setminus \{a_j\}$ for some $j\geq 8$. So, using \eqref{39} we get 
$$
|f'(z)|\geq \frac{b_j}{|z-a_j|^2} -\!\! \sum_{k\neq j,\, k\geq 8}\frac{b_k}{|z-a_k|^2} 
\geq \frac{b_j}{b_j^2} - 9\sum_{ k\geq 8}k^2 e^{-k}\geq e^8-1 >0.
$$
Therefore we have again $\cri\subset \cU$ and $f(\cri) \subset \ov \D (0, 2)$.
Arguing exactly as in the proof of Theorem~\ref{prop ex1}, it also follows that all possible asymptotic values of $f$ belong to $\ov \D (0, 2)$. This shows that $f\in \cB$ and that $\infty$ is not an asymptotic value of $f$.

Thus, the function $f$ is in the class of functions considered in Theorem \ref{theo main}, and this theorem implies that $\HD(\escf) =\d $. It remains to determine the value of $\d$. Since all the poles of $f$ are simple, the series determining this value is
$$
\sum_{j\geq 8} \Big(\frac{b_j}{a_j^{2}}\Big)^t =\sum_{j\geq 8} \Big(\frac{e^{-j}}{(\log j)^{2}}\Big)^t.
$$
It converges for every $t>0$. In consequence (using Theorem \ref{theo main}), $\HD(\escf) =\d =0$. 
\epf


\bibliographystyle{plain}

\end{document}